\documentclass[11pt]{article}
\usepackage{graphicx}
\usepackage{amsthm}
\newtheorem{theorem}{Theorem}

\newtheorem{corollary}[theorem]{Corollary}

\newcommand{\p}{\mbox{${\cal P}_n$}}
\newcommand{\fs}{\mbox{$\mathcal{S}_n$}}
\newcommand{\comp}{\mbox{$\mathcal{C}_n$}}
\newcommand{\vn}{\mbox{$\mathcal{V}_n$}}
\newcommand{\vnk}{\mbox{$\mathcal{V}_{n,k}$}}
\begin{document}
\title{Valleyless Sequences}
\author{Robert G. Rieper\\William Paterson University  \and Melkamu Zeleke
\thanks{Corresponding author.} \\William Paterson University}
\date{}
\maketitle
\begin{abstract}
Valleyless sequences of finite length $n$ and maximum entry $k$
occur in tree enumeration problems and provide an interesting
correspondence between permutations and compositions. In this
paper we introduce the notion of \emph {valleyless} sequences,
explore the correspondence and enumerate them using the method of
generating functions.
\end{abstract}

\section{Introduction}

Throughout this article we use the notation $\fs$ to denote finite
sequences of length $n$. A length $n$ sequence $s=s_1,s_2,\ldots
,s_n$ of positive integers is a permutation whenever each $s_i$ is
a distinct member of the set $\{1,2,\ldots,n\}$.  In this case we
write the sequence $s$ as a word $\pi=\pi_1\pi_2\ldots\pi_n$ and
denote the set of all permutations of length $n$ as $\p$ . In the
general case, the sequence of positive integers $s$ can be
associated with a composition of the integer $s_1+s_2+\ldots+s_n$.
The composition has $n$ parts where $s_i$ is the size of the
$i^{th}$ part. We let $\comp$ be the set of all compositions of
the positive integer $n$.

A sequence is defined herein to be \emph{valleyless} provided
$1\le i<j<k\le n$ implies $s_j\ge \min{\{s_i,s_k\}}$.  A graph of
the pairs $(i,s_i)$ reveals that a valleyless sequence never has a
valley.  As a sequence, it is either nondecreasing, nonincreasing,
or is nondecreasing to a point and thereafter nonincreasing.  We
let $\vn$ be the set of all length $n$ valleyless sequences of
positive integers and $\vnk$ the set of all length $n$ valleyless
sequences of positive integers with maximum part $k$.

Valleyless sequences occur in tree enumeration
problems~\cite{RZ99} and provide an interesting correspondence
between permutations and compositions.  We explore this
correspondence first and then enumerate valleyless sequences in
three different ways.

\section{Valleyless permutations}
A permutation $\pi=\pi_1\pi_2\ldots\pi_n\in \p$ is uniquely
determined by its inversion table $I(\pi)=(a_1,a_2,\ldots,a_n)$.
The entry $a_k$ is the number of symbols in the word $\pi$ to the
left of $k$ that are greater than $k$. Thus, for all $k$, $0\le
a_k\le n-k$.  The total number of inversions of $\pi$ is
$i(\pi)=a_1+a_2+\ldots +a_n$ and $$ \sum_{\pi \in
\mathcal{S}_n}q^{i(\pi)}=(1+q)(1+q+q^2)\cdots(1+q+q^2+\ldots
+q^{n-1}).$$  As an example we have $\pi=2731546$, $a_1=3$,
$I(\pi)=(3,0,1,2,1,1,0)$, and $i(\pi)=8$.

A finite sequence $s=s_1,s_2,\ldots ,s_n$ viewed as a composition
$s_1+s_2+\ldots+s_n=N$ can be encoded as a unique subset of
$\{1,2,\ldots,N-1\}$ by
$\Theta(s)=\{s_1,s_1+s_2,\ldots,s_1+s_2+\ldots +s_{n-1}\}$.  For
example, the sequence $s=2,3,2,1,1$ is the composition
$2+3+2+1+1=9$ with $\Theta(s)=\{2,5,7,8\}$.  The process is
reversible so that there are $2^{N-1}$ different compositions of
$N$.

  Next, we present a result
relating \emph{valleyless} permutations and compositions.

\begin{theorem}
\label{t:one} The subset of valleyless permutations in $\p$ is in
one-to-one correspondence with the compositions of $n$. That is,
$|\vn\cap \p |=|\comp|$.
\end{theorem}

\begin{proof}Let $\pi=\pi_1\pi_2\ldots\pi_n$ be valleyless, then either
$\pi_1=1$ or $\pi_n=1$.  Add one to every entry of $\pi$ and then
append 1 to the beginning or end.  The new permutation is
valleyless of length $n+1$.  The process is reversible so there
are twice as many valleyless permutations of length $n+1$ as there
are of length $n$.  Since $|\mathcal{P}_1\cap \mathcal{S}_1 |=1$
we have $|\vn\cap \p |=2^{n-1}=|\comp|$ and we are done.
\end{proof}
The statement of the theorem implies a stronger correspondence
than the cardinality result.  There is indeed such a
correspondence which we present now.  Our goal in presenting the
above proof was to introduce the reader to the technique of
altering a sequence by adding one to every entry and then
appending a 1 at the beginning or end.  We will use this technique
in the future.  The stronger result is obtained from the
following.

\begin{theorem}
\label{t:two} A permutation $\pi=\pi_1\pi_2\ldots\pi_n$ with
inversion table $I(\pi)=(a_1,a_2,\ldots,a_n)$ is valleyless if and
only if $a_k=0$ or $n-k$, $k=1,2,\ldots ,n$.

\end{theorem}
\begin{proof}
Suppose that a permutation $\pi=\pi_1\pi_2\ldots\pi_n$ with
inversion table $(a_1,a_2,\ldots,a_n)$ is valleyless. If $a_k>0$
for some $k$, $1\le k\le n$, then there is at least one symbol in
$\pi$ to the left of $k$ that is greater than $k$. If there is any
symbol in $\pi$ greater than $k$ to the right of $k$, then $k$
would be a valley contradicting our assumption. Hence all numbers
greater than $k$ must be to the left of $k$ and $a_k=n-k$.
Conversely, suppose that $I(\pi)=(a_1,a_2,\ldots,a_n)$ and $a_k=0$
or $n-k$, $k=1,2,\ldots ,n$. We want to show that $\pi$ belongs to
$\vn\cap \p$. Suppose not, then there is at least one symbol
${\pi}_k$ in $\pi$ such that
\[
{\pi}_k < {\pi}_i,\mbox{ for some } i<k \mbox{ and } {\pi}_k <
{\pi}_j,\mbox{ for some } j>k.
\]
Then $a_k\ne 0$ and $a_k < n-k$ contradicting our assumption.
Hence $\pi$ belongs to $\vn\cap \p$ and we are done.
\end{proof}
Note that there are two choices for each $a_k$ except for $k=n$ in
which case there is only one.  Thus, there are $2^{n-1}$
valleyless permutations of length $n$.  If we again use $q$ to
record the number of inversions in a permutation, then we have
\begin{corollary}
\label{c:one} $$\sum_{\pi \in \vn \cap \p}
q^{i(\pi)}=(1+q)(1+q^2)\cdots (1+q^{n-1}).$$
\end{corollary}
\section{Permutations with a valley}
One of the fundamental statistics associated with a permutation
$\pi =\pi _1\cdots \pi_n$ is the descent set $D( \pi ) = \{i|\pi_i
> \pi_{i+1}\}.$ See ~\cite{RS86} for the enumeration of permutations by their descent set.
We say that a length $n$ sequence $s=s_1,s_2,\cdots,s_n$ has
exactly $k$ valleys if $$|\{j: s_j <
min\{s_{j-1},s_{j+1}\}\}|=k.$$ In this section we obtain a
recursive formula for the generating function of permutations with
exactly $k$ valleys which we denote by $g_k(x)$ and use the
formula to reproduce the table of permutations with peaks in
~\cite{DKB66}.
\\

\noindent We have shown in the previous section that $g_0(x)$
satisfies the relation
\begin{eqnarray}
 g_0(x)=\frac x{1-2x}.
\end{eqnarray}
\noindent One can easily see that length $n+1$ permutations with
exactly $k+1$ valleys can be obtained from:
\begin{enumerate}
\item length $n$ permutations with exactly $k+1$ valleys by adding
$1$ to every entry and then appending a $1$ either at the
beginning, end, to the left or to the right of the $k+1$ valleys.
\item length $n$ permutations with exactly $k$ valleys by adding $1$ to every entry and
then inserting a $1$ in one of the $(n-1)-2k = n-(2k+1)$ middle
positions.
\end{enumerate}
\noindent Hence
\begin{eqnarray}
g_{k+1}(x)=(2(k+1)+2)xg_{k+1}(x)+x^{2k+3}D_x\left (\frac
{g_k(x)}{x^{2k+1}}\right ), \mbox{ for } k\ge 0.
\end{eqnarray}
\noindent Substituting $f_k(x)=\frac{g_k(x)}{x^{2k+1}}$ in $(2)$
we obtain a system of linear differential equations
\begin{eqnarray}
f_{k+1}(x)=\frac 1{1-(2(k+1)+2)x} D_x(f_k(x))
\end{eqnarray}
which can also be written in matrix form as
\begin{eqnarray*}
\left [ \begin{array}{c} f_0 \\ f_1 \\ f_2 \\ f_3 \\ f_4 \\ \vdots
\end{array} \right ] = \left [ \frac 1{1-2x},  \frac 1{1-4x},  \frac
1{1-6x}, \frac 1{1-8x},\ldots \right ]\cdot D_x\left [
\begin{array}{c} x \\ f_0 \\ f_1 \\ f_2 \\ f_3 \\ \vdots
\end{array} \right ]
\end{eqnarray*}
 \noindent Even though $(3)$ is an elegant recurrence relation
 just like many known simple relations for the number of
 partitions of an integer, we do not know how to solve it at this time.
\begin{center}
\begin{tabular} {|c|p{10cm}|} \hline
k & $g_k(x)$ \\ \hline 0 &  $\frac {x}{1-2x}$
\\[2ex] \hline 1 &  $\frac {2x^3}{(1-2x)^2(1-4x)}$   \\[2ex] \hline 2 &  $\frac
{16x^5(1-3x)}{(1-2x)^3(1-4x)^2(1-6x)}$
\\[2ex] \hline 3 & $\frac
{16x^7(17-184x+636x^2-720x^3)}{(1-2x)^4(1-4x)^3(1-6x)^2(1-8x)}$
\\[2ex] \hline
4 & $\frac
{256x^9(31-788x+8096x^2-43132x^3+126072x^4-192672x^5+120960x^6)}{(1-2x)^5(1-4x)^4(1-6x)^3(1-8x)^2(1-10x)}$
\\[2ex] \hline
\end{tabular}\\[.75ex]
{\it Table 1: Generating functions of permutations with valleys}
\end{center}

\noindent The formal Taylor series expansion of the functions in
in Table 1  reproduces the table of permutations with peaks in
~\cite{DKB66}.
\begin{center}
\begin{tabular} {|r||r|r|r|r|r|r|r|r|r|r|r|r|r|r|} \hline
k/n & 1 & 2 & 3 & 4 & 5 & 6 & 7 & 8 & 9 & 10
\\ \hline\hline 0 & 1 & 2 & 4 & 8 & 16& 32 & 64 & 128 & 256 & 512
\\ \hline 1 & & & 2 & 16 & 88 & 416 & 1824 & 7680 & 31616 & 128512 \\ \hline
2 & & & & & 16 & 272 & 2880 & 24576 & 185856 & 1304832
\\ \hline 3 & & & & & & & 272 & 7936 & 137216 & 1841152
\\ \hline
4 & & & & & & & & & 7936 & 353792  \\ \hline
\end{tabular}\\[.75ex]
{\it Table 2: The number of permutations of length $n$ with
exactly $k$ valleys}
\end{center}
If we denote the number of permutations of $n$ numbers with $k$
valleys by $P(n,k)$, the functional recursion (2) can also be
written as
\begin{eqnarray}
P(n,k) = 2(k+1)P(n-1,k) + (n-2k)P(n-1,k-1)
\end{eqnarray}
and notice the remarkable similarity of the recurrence relation
(4) with that of Eulerian numbers ~\cite{GKP89,WS96}
$$E(n,k)=(k+1)E(n-1,k)+(n-k)E(n-1,k-1).$$
\section{Valleyless sequences}
In this section we enumerate valleyless sequences of length $n$
and maximum entry $k$ which we have denoted by $\vnk$ using the
method of generating functions. Table 3 shows the cardinalities of
$\vnk$ for $1\le n\le 6$ and $1\le k\le 5$.

\begin{center}
\begin{tabular} {|r||r|r|r|r|r|} \hline
n/k & 1 & 2 & 3 & 4 & 5 \\ \hline\hline 1 & 1 & 1 & 1 & 1 & 1
\\ \hline 2 & 1 & 3 & 5 & 7 & 9  \\ \hline 3 & 1 & 6 & 15 & 28 &
45
\\ \hline 4 & 1 & 10 & 35 & 84 & 165  \\ \hline 5 & 1 & 15 & 70 & 210 &
495
\\ \hline 6 & 1 & 21 & 126 & 462 & 1287 \\ \hline
\end{tabular}\\[.75ex]
{\it Table 3: Valleyless Sequences of Length $n$ with maximum
entry $k$ }
\end{center}

Let $V(x,y)$ be a generating function which enumerates valleyless
sequences of length $n$ and maximum entry $k$. Then
\[V(x,y)=\sum_{n,k} V_{n,k} x^ny^k = 1+xy+xy^2+\ldots, \]
where $V_{n,k} = |\vnk|$. As pointed out in the proof of Theorem
1, any valleyless sequence of length $n$ and maximum entry $k$ can
be obtained from the base case:
\[1, 11, 111, 1111, \ldots \]
by adding one to every entry and then appending a $1$ at the
beginning or end. Hence $V(x,y)$ satisfies the recursion
\begin{eqnarray}
V(x,y)& = & \underbrace{\{x+x^2y+x^3y+\ldots\}}_{\mbox{base case}}
+\overbrace{\{yV(x,y)\}}^{\mbox{adding 1}}\times
\underbrace{\{1+2x+3x^2+\ldots\}}_{\mbox{appending a 1}} \nonumber
\\ & = & \frac {xy}{1-x} + yV(x,y)\frac{1}{(1-x)^2}. \nonumber
\end{eqnarray}
Therefore,
\begin{eqnarray} V(x,y)\left(1-\frac{y}{(1-x)^2}\right )= \frac {xy}{1-x} \mbox{
and } V(x,y)=\frac {xy}{1-x-\frac {y}{1-x}}. \end{eqnarray} On the
other hand, it is well known ~\cite{AT95} that
\[\sum_{k=0}^{\infty} {n-1+k \choose k} y^k = \frac
{1}{(1-y)^n}.\] Using this closed form and the 'snake oil'
~\cite{HW90} method, we observe that
\begin{eqnarray}
& & \sum_{n=1}^{\infty}\left (\sum_{k=1}^{\infty}{n-1+2(k-1)
\choose 2(k-1)} y^k\right )x^n \\ & = & \sum_{n=1}^{\infty} \frac
y2 \left ( \frac {1}{(1-\sqrt{y})^n} +  \frac {1}{(1+\sqrt{y})^n}
\right ) x^n \nonumber \\ & = &  \frac y2
\sum_{n=1}^{\infty}{\left ( \frac {x}{1-\sqrt{y}} \right )}^n +
\frac y2 \sum_{n=1}^{\infty}{\left ( \frac {x}{1+\sqrt{y}} \right
)}^n \nonumber \\ & = & \frac y2 \left ( \frac {x}{1-\sqrt{y}-x} +
\frac {x}{1+\sqrt{y}-x}\right )\nonumber \\ & = & \frac
{xy}{(1-x)-\frac {y}{1-x}}.
\end{eqnarray}
Hence from $(5)$ and $(6)$ one can see that
\begin{theorem}
\label{t:3} The number of valleyless sequences of length $n$ with
maximum entry $k$ is \[ V_{n,k} = {n-1+2(k-1)\choose 2(k-1)}. \]
\end{theorem}
\section{q-Analog}
If we use $q$ to record the sum of the entries of members of
$\vnk$ and let $V(x,q,y)$ its generating function, then
\[V(x,q,y)=\sum_{n,p,k}  V_{n,p,k} x^nq^py^k. \]
Using the same argument as in the previous case, we see that
$V(x,q,y)$ satisfies the functional recursion
\begin{eqnarray}
V(x,q,y) & = & \frac {xqy}{1-xq} + yV(xq,q,y)\frac {1}{(1-xq)^2}
\nonumber \\ \Leftrightarrow (1-xq)^2V(x,q,y) & = &
xqy-x^2q^2y+yV(xq,q,y).
\end{eqnarray}
Equation $(7)$ is a linear $q-$difference equation and we seek
series solution of the form
\begin{eqnarray}
V(x,q,y)= \sum_{n=0}^{\infty} b_n(x,q)y^n.\end{eqnarray}
Substituting $(8)$ into $(7)$ and comparing coefficients of $y^n$
we obtain:
\begin{eqnarray}
b_n(x,q)=\frac {b_{n-1}(xq,q)}{(1-xq)^2}\mbox{ for } n>1,
\end{eqnarray}
where $ b_0=0 \mbox{ and } b_1=\frac {xq}{1-xq}.$\\ Repeated
application of equation (9) gives the explicit solution
\begin{eqnarray}
b_n=\frac {xq^n(1-xq^n)}{(xq)^{2}_{n}},
\end{eqnarray}
where the common $q-$symbol \[(xq)_{n} = \prod_{i=1}^{n}
(1-xq^i).\] Therefore,
\begin{eqnarray}
V(x,q,y)& = &\frac {xq}{1-xq} y + \sum_{n=2}^{\infty} \frac
{xq^n(1-xq^n)}{(xq)_n^2}y^n \nonumber \\ & = & \sum_{n=1}^{\infty}
\frac {xq^n(1-xq^n)}{(xq)_n^2}y^n
\end{eqnarray}
An implementation of $(11)$ using Maple\footnote{Send e-mail to
{\bf zelekem@wpunj.edu} to obtain the Maple Package {\bf bxq}
accompanying this paper.} shows that the number of valleyless
sequences of length $10$, sum of entries $20$, and maximum part
$5$, for instance, is 325.

\section{A nonlinear relation} Expanding $V(x,q,y)$ in $x$ as $V(x,q,y)=\sum a_n(q,y)x^n$,
substituting it into $(5)$, and comparing coefficients of $x^n$,
we obtain a three term recurrence relation
\begin{eqnarray}
 a_n=\frac {2qa_{n-1}-q^2a_{n-2}}{1-yq^n} \mbox{ for } n>2
\end{eqnarray}
where $a_0=0, a_1=\frac {yq}{1-yq},\mbox{ and } a_2=\frac
{q^2y+q^3y^2}{{(yq;q)}_2}$.

This nonlinear recurrence relation provides yet another way of
enumerating valleyless sequences. Indeed, one can obtain
valleyless sequences of length $n$ from sequences of length $n-1$
by adding a 1 at the beginning or end, and then appending 1s to
every entry. This process contributes $\frac {2q}{1-yq^n}
a_{n-1}(q,y)$ to $a_n(q,y)$. However, sequences with the same
first and last entry, for example $23522$, are counted twice in
this algorithm and we need to subtract $\frac {q^2}{1-yq^n}
a_{n-2}(q,y)$, the number of valleyless sequences of length $n$
with the same first and last entry. Hence,
\[a_n(q,y) = \frac {2q}{1-yq^n}
a_{n-1}(q,y) - \frac {q^2}{1-yq^n} a_{n-2}(q,y). \] Unlike (9)
equation (12) is nonlinear and we do not know its explicit
solution.

\section{Conclusion} We came across valleyless sequences while
enumerating ordered trees by their number of leaves, total path
length and number of vertices. Given a positive integer $n$,
consider its composition or ordered partition and draw a
composition tree corresponding to $n$. Then try to construct all
possible ordered trees of total path length $n$ by identifying
non-terminal nodes or vertices of the composition tree.

A sequence obtained from an ordered partition of a positive
integer $n$ by taking the minimum of two consecutive entries gives
a sequence without isolated valley and we hope that the ideas
developed in this article will provide an alternative means to
solve the ordered tree enumeration problem.

\section{Acknowledgements}
The second author acknowledges research release time from William
Paterson University and summer support from The Center for
Research of the College of Science and Health.
\bibliographystyle{amsplain}
\bibliography{valleybib}

\vspace{.5in} \noindent\texttt{jrieper@cybernex.net}\\
\texttt{zelekem@wpunj.edu}\\
\end{document}